\documentclass[12pt,leqno]{amsart}
\usepackage{amssymb, graphicx, amscd, amsmath, amsthm, 
array, fancyvrb, relsize, paralist,url}
\usepackage{fullpage}
\usepackage{colonequals}
\usepackage[all]{xy}

\makeatletter
\g@addto@macro\bfseries{\boldmath}
\makeatother

\theoremstyle{plain}
\newtheorem{thm}{Theorem}[section]
\newtheorem{lem}[thm]{Lemma}

\newtheorem{exa}[thm]{Example}

\newtheorem*{thm*}{Theorem}

\def\C{\mathcal{C}}
\def\S{\mathcal{S}}
\def\Gal{\mathrm{Gal}}

\def\bcases{\begin{cases}}
\def\ecases{\end{cases}}

\newcommand{\jac}[2]{\left(\displaystyle{\frac{#1}{#2}}\right)}

\newcommand{\s}{\sigma}
\newcommand{\Q}{\mathbb{Q}}
\newcommand{\Z}{\mathbb{Z}}
\newcommand{\OO}{\mathcal{O}}

\newcommand{\PP}{\mathfrak{p}}

\renewcommand{\a}{\alpha}
\renewcommand{\b}{\beta}
\renewcommand{\c}{\gamma}

\newcommand{\w}{\omega}

\newcommand{\bp}{\overline{\pi}}
\newcommand{\bE}{\overline{E}}
\newcommand{\p}{\pi}

\numberwithin{equation}{section}

\begin{document}

\title{Cubic residuacity of real quadratic integers}

\author{Ron Evans}
\address{Department of Mathematics, UCSD\\
La Jolla, CA  92093-0112}
\email{revans@ucsd.edu}
\urladdr{https://mathweb.ucsd.edu/~revans}

\author{Mark Van Veen}
\address{2138 Edinburg Avenue\\
Cardiff by the Sea, CA 92007}
\email{mavanveen@ucsd.edu}
\urladdr{}



\date{Nov 2025}

\begin{abstract}
Given a real quadratic integer $u=A+B\sqrt{D}$ with cubic norm,
we identify all the classes in a related form class group that
represent primes $p$ for which
$u$ is a cubic residue mod $p$.  A special case of this result was conjectured
in a 2025 paper of Evans, Lemmermeyer, Sun, and Van Veen.
\end{abstract}

\maketitle

\section{Introduction}
For squarefree $D>1$, write $\OO$ for the ring of integers in
$\Q(\sqrt{D})$.
Let $\S_D$ denote the set 
consisting of those $u = A+B\sqrt{D}$ with $A,B \in \Z$
such that the norm $N(u)$ is a cube, $u$ itself is
not a cube in $\OO$, and $(A,B)$ is not divisible by the cube
of a prime. 
A closely related set $S_D$ 
was introduced in \cite[Section 7]{ELSV}.
The difference between the two sets $S_D$ and $\S_D$ is that
those elements $\mu$ in $S_D$ of the form $(x+y\sqrt{D})/2$
with odd $x,y$ have been converted to $8\mu$ in $\S_D$
to ensure integer coefficients without affecting cubic residuacity.
Some examples of elements in $\S_D$  for $D \in \{2,3,5,6,7\}$ 
are $19+3\sqrt{2}$, $1342+99\sqrt{3}$, $3047+176\sqrt{5}$,
$1633+437\sqrt{6}$, $232+319\sqrt{7}$, with norms $7^3$, $121^3$,
$209^3$, $115^3$, $-87^3$, respectively.

For $u \in \S_D$, define $C = C(u)$ to be the
product of the distinct odd primes that 
divide $(A,B)$ but not $D$.  For example,  $C(u)=5$ for $u=110+10\sqrt{41}$,
while $C(u)=1$ for $u=25+5\sqrt{30}$.

Corresponding to $u \in \S_D$, let $d = -27C^2D$, where $C=C(u)$.
Let $H = H(4d)$ be the form class group of discriminant $4d$   
consisting of the
classes of primitive binary quadratic forms $[a, 2b,c]$ with $a>0$
and $d=b^2 -ac$.

The primary purpose of this paper is to prove Theorem 1.1 below,
which identifies the
classes in $H$
that represent primes $p \equiv 1 \pmod{3}$  for which
$u=A+B\sqrt{D}$  is a cubic residue mod $p$.
Note that $\sqrt{D}$ mod $p$ exists for such $p$,
since $d$ is a square mod $p$ \cite[Prop. 12.3.1]{G}, \cite[Sec. 6.2]{HK}.
The cubic residuacity of $u$ does not depend on the choice of
$\sqrt{D}$ mod $p$, since $N(u)$ is a cube.
Theorem 1.1 is proved in Section 3.

Versions of Theorem 1.1 for fundamental units $u$ 
are due to Sun \cite[Section 5]{Sun1}.
The special case of our Theorem 1.1 applied to the principal class
in $H$ proves the conjecture in \cite[Conj. 7.8]{ELSV}.
That conjecture asserted that the real cube root
of every $u \in \S_D$ lies in the ring class field for the order
$\Z[3C\sqrt{-3D}]$ in $\Q(\sqrt{-3D})$. This means that
for primes $p$ of the form $p=x^2+27C^2Dy^2$, \emph{every} $u \in \S_D$ 
with $C(u)=C$
is a cubic residue mod $p$.
We remark that the same is not true for primes
$p$ of the forms $p=x^2 +3C^2Dy^2$, $p=x^2 +27CDy^2$, or $p=x^2+27Dy^2$.
For example, $17+51\sqrt{2} \in \S_2$ is not cubic mod any of the primes
$p=1759 = 25+3C^2D$, $p=919=1+27CD$, $p=79=25+27D$.

Define $\w = (-1+i\sqrt{3})/2$.
Our proof of Theorem 1.1 is based on a modification of an ingenious
argument of Sun \cite{Sun2}, which invokes
properties of the cubic residue symbol
$\jac{*}{*}_3$
for the Eisenstein integers $\Z[\w]$ \cite[Ch. 8]{BEW}, \cite[Sec. 7.1]{HK}.

\begin{thm}\label{Theorem 1.1}
Let $u =A+B\sqrt{D}\in \S_D$ and 
let $p \equiv 1 \pmod{3}$ be a prime represented
by a class $[a,2b,c] \in H$, so that $p=ax^2 +2bxy+cy^2$ for some integers
$x,y$.  Assume that
\begin{equation}\label{eq 1.1}
(ap, 3N(u)Dy)=1 .
\end{equation}
Write $C=C(u)$. Then $u$ is a cubic residue mod $p$ if and only if
\begin{equation}\label{eq 1.2}
\jac{9A -(B/C)b(1+2\w)}{a}_3 =1.
\end{equation}
\end{thm}

\noindent
Observe that the symbol in (1.2) is independent of the choice of 
prime $p$ represented by $[a,2b,c]$.

After the following comments
demonstrating the ``insignificance" of the restriction (1.1),
we provide an example for Theorem 1.1.
Let $\C$ be any class in $H$.  By \cite[pp. 149--150]{Cox},
$\C$ represents infinitely many primes.  Suppose that $a$ is one of these
primes with $a > 3|N(u)|D$.  Then by \cite[eq. 12.1.4]{G}, \cite[Sec. 6.2]{HK},
$\C$ can be written as
$[a, 2b,c]$ for some integers $b,c$. Let $p>a$  be a prime
with $p=ax^2 + 2bxy +cy^2$.   Since
$(ax+by)^2 -dy^2=ap$,
we have $(ap,y)=1$, so that $(ap, 3N(u)Dy)=1$.
In short, there always exists a form $(a,2b,c) \in \C$
that satisfies (1.1) for infinitely many primes $p$ represented by $\C$.

\begin{exa}\label{ex 1.2}
Let $u=6+3\sqrt{7}$, so that $N(u)=-27$, $D=7$, and $C=3$.
The class group $H$ consists of $36$ classes of discriminant
$4d = -6804$,  half of which represent primes $p \equiv 1 \pmod{3}$.
Here (1.2) 
holds for six classes in $H$, namely
\begin{equation}\label{eq 1.3}
[1,0,1701], \ [7,0,243],\ [19,\pm 6,90], \ [22,\pm 18,81].
\end{equation}
The classes in (1.3) are precisely those in $H$ 
that represent primes $p \equiv 1 \pmod{3}$
for which $u$ is a cubic residue mod $p$.
As is explained in Section 4,
it should come as no
surprise that the classes in (1.3) form a subgroup of index $6$
in the class group, 
since the same is true in general for the classes satisfying (1.2).
\end{exa}

In the sequel, let
$u = A +B \sqrt{D} \in \S_D$, $[a,2b,c] \in H$ and
$f=f(x,y) = ax^2 +2bxy+cy^2$ for integers $x,y$
(where $f$ need not be prime).
We have
\begin{equation} \label{eq 1.4}
(ax+by)^2 -dy^2=af(x,y).
\end{equation}
When $3 \nmid f(x,y)$, define
\begin{equation}\label{eq 1.5}
L=L(u)=L(u,f)=\jac{9Ay+(B/C)(ax+by)(1+2\w)}{f(x,y)}_3,
\end{equation}
and
when $3 \nmid a$, define
\begin{equation}\label{eq 1.6}
R=R(u)=R(u,f)=\jac{9A-(B/C)b(1+2\w)}{a}_3.
\end{equation}
In  Theorem 2.4, we prove that $L=R$ whenever
\begin{equation}\label{eq 1.7}
(af(x,y), 3N(u)Dy) = 1.
\end{equation}
This result is crucial for the proof of Theorem 1.1.

\section{$L=R$}

Corresponding to $u = A + B\sqrt{D} \in \S_D$, 
define $Q: = (A, B/C)$.  Assume that (1.7) holds.
The object of this section is to prove that $L=R$
(see Theorem 2.4).
Our proof makes use of the following three well-known properties
of the cubic residue symbol for integers $m, m', n, n'$
\cite[Sec. 7.1]{HK}.
When $3 \nmid m$,
\[
\jac{1+2\w}{m+3n\w}_3 = \w^{jn}, \quad \mbox{where} \ j = \jac{m}{3} .
\]
When $(n,3m)=1$,
\[
\jac{m}{n}_3 = 1.
\]
When $3 \nmid m m'$,
\[
\jac{m+3n\w}{m' +3n'\w}_3 = \jac{m'+3n'\w}{m +3n\w}_3.
\]

We proceed  with three lemmas.

\begin{lem}\label{Lemma 2.1}
Let $q$ be a prime dividing $Q$. Write
$s=v_q(A)$ and  $r=v_q(B/C)$,  where $v_q(n)$
denotes the largest power of $q$ dividing $n$.
If either $q>2$ or $q=2$ with $2\mid D$, then $s>r$.
If $q=2$ with $2 \nmid D$, then $s=r$.
\end{lem}

\begin{proof}
By definition of $q$, the integers $r,s$ are positive.
Recall that $N(u) = A^2 -B^2D$.

First suppose that $q >2$.
If $q\mid D$ (so that $q\nmid C$), then $v_q(B^2D)=2r+1$.
If $q\nmid D$ (so that $q\mid C$), then $v_q(B^2D)=2r+2$.
Suppose for the purpose of contradiction that $s \le r$.
Then $v_q(N(u)) =2s$.   Since $N(u)$ is a cube, $3 \mid s$,
so that $r \ge s \ge 3$.  Then $(A,B)$ is divisible by $q^3$,
which contradicts the definition of the set $\S_D$.
Therefore $s > r$.

Next suppose that $q=2$, so that $3 \mid v_2(N(u))$.
If $D$ is even, then arguing by contradiction  as above, we see that $s>r$.
If $D$ is odd, then clearly we must have $r=s \le 2$ (with strict inequality
when $D \equiv 3 \pmod 4$).

\end{proof}

Define $u_2 = (A_2 + B_2\sqrt{D})/T$, where $A_2 = BD^2$ , $B_2 = AD$,
and $T$ is the largest integer cube dividing $(A_2, B_2)$.
We have $u_2 = (\sqrt{D})^3 u/T$ and $N(u_2) = -D^3 N(u) /T^2$, so
$u_2 \in \S_D$.  Note that $C=C(u) = C(u_2)$.
By (1.7), $(T, af(x,y))=1$, because any prime divisor of $T$
must also divide $D$.

\begin{lem}\label{Lemma 2.2}
$R(u) = R(u_2)$.
\end{lem}

\begin{proof}
Since $b^2-ac=d=-27C^2D$, we have $b \equiv  3C\s \pmod{a}$,
where $\s$ is a square root of $-3D$ mod $a$.  Then
\[
R(u_2) = \jac{9BD -(A/C)b(1+2\w)}{a}_3
=\jac{9BD -3A\s(1+2\w)}{a}_3,
\]
so that
\[
R(u_2)=\jac{\s}{a}_3 \jac{-3B\s-3A(1+2\w)}{a}_3.
\]
Since $\jac{\s}{a}_3=1$ and $\jac{1+2\w}{a}_3=1$, we obtain
\[
R(u_2)=\jac{9A-3B\s(1+2\w)}{a}_3=\jac{9A-(B/C)b(1+2\w)}{a}_3=R(u).
\]
\end{proof}

\begin{lem}\label{Lemma 2.3}
$L(u) = L(u_2)$.
\end{lem}

\begin{proof}
By (1.4),
we have $ax+by \equiv  3yC\tau \pmod{f(x,y)}$,
where $\tau$ is a square root of $-3D$ mod $f(x,y)$.  Then
\[
L(u_2) = \jac{9BDy +(A/C)(ax+by)(1+2\w)}{f(x,y)}_3
=\jac{9BDy +3Ay\tau(1+2\w)}{f(x,y)}_3,
\]
so that
\[
L(u_2)=\jac{\tau}{f(x,y)}_3 \jac{-3By\tau+3Ay(1+2\w)}{f(x,y)}_3.
\]
Since $\jac{\tau}{f(x,y)}_3=1$ and $\jac{(1+2\w)}{f(x,y)}_3=1$, we obtain
\begin{equation}\label{eq 2.1}
L(u_2)=\jac{9Ay+3yB\tau(1+2\w)}{f(x,y)}_3=
\jac{9Ay+(B/C)(ax+by)(1+2\w)}{f(x,y)}_3=L(u).
\end{equation}
\end{proof}

In preparation for the proof of Theorem 2.4,
we'll need the following additional facts.
First we show that $L$ and $R$ are nonzero.
By (1.4), the norm of the numerator of $L$ in (1.5) is
congruent modulo $af(x,y)$ to $81N(u)y^2$, which by (1.7)
is coprime with $af(x,y)$.  Thus $L$ is nonzero.
Since $b^2 = ac+d$ is congruent modulo $a$ to $d$,
the norm of the numerator of $R$ in (1.6) is congruent mod $a$
to $81N(u)$, which is coprime with $a$.   Thus $R$ is nonzero.

Recall that $Q = (A, B/C)$.  By (1.7), $(Q, af(x,y))=1$,
since $Q^2 \mid N(u)$.  We proceed to show that
\begin{equation}\label{eq 2.2}
(ax+by, 3yQ') = 1,
\end{equation}
where $Q'$ denotes the odd part of $Q$.
If some prime were to divide $(ax+by, 3y)$,
then it would also divide $af(x,y)$ by (1.4),
contradicting (1.7).  Next let $q$ be a prime
dividing $Q'$. If $q \mid D$, then $q \mid d$.  If $q \nmid D$,
then $q \mid C$, so $q \mid d$ in either case.  Therefore $q \nmid (ax+by)$,
otherwise $q$ would divide $af(x,y)$ by (1.4), contradicting the fact that
$(Q, af(x,y))=1$. Thus (2.2) holds.

\begin{thm}
Suppose that (1.7) holds. Then  $L=R$.
\end{thm}

\begin{proof}
Let $\a=A/Q$ and $\b=B/(CQ)$, so that $(\a,\b)=1$.
Let $\c=(y, \beta)$, $\b_1=\b/\c$, $y_1=y/\c$, so that $(\b_1, y_1)=1$.
By (1.7),  $(\c,af(x,y)=1$. Define $M$ by
\begin{equation}\label{eq 2.3}
M =\b_1(ax+by)-3y_1\a(1+2\w).
\end{equation}

We now show that $(N(M), af(x,y))=1$.  Taking the norm of $\c M$, we have
\begin{equation}\label{eq 2.4}
\c^2 N(M) = \b^2 (ax+by)^2 +27 y^2 \a^2, 
\end{equation}
so by (1.4),
$\c^2 N(M) = \b^2 af(x,y) +27 y^2 (\a^2-\b^2C^2D)$.
Therefore $N(M) \equiv  27 y_1^2 N(u)/Q^2 \pmod{af(x,y)}$,
so by (1.7), $(N(M),af(x,y))=1$.

We next show that 
\begin{equation}\label{eq. 2.5}
(N(M),Q)=1.
\end{equation}
Let $q$ be a prime dividing $Q$.  If $q$ is odd, then
by Lemma 2.1, $q \nmid \b$ and $q \mid \a$, so
by (2.2) and (2.4),  $(N(M), Q')=1$.
Now suppose that $q=2$.  
To prove (2.5), it remains to show that $N(M)$ is odd.
Since $Q$ is even, $af(x,y)$ must be
odd by (1.7).  First suppose that $D$ is even. Then $(ax+by)$
is odd by (1.4).  By Lemma 2.1, 
$2 \nmid \b$ and $2 \mid \a$, so $N(M)$ is odd by (2.4).
Next suppose that $D$ is odd.  Then $d$ is also odd,
since $C$ is odd. By Lemma 2.1,  both $\a$ and $\b$
are odd.  Therefore by (1.4), $(ax+by)$ and $y$ have opposite parity.
It now follows from (2.4) that $N(M)$ is odd.
This completes the proof of (2.5).

Since $(Qy, f(x,y))=1$ by (1.7), it follows from (1.5) that
\[
L=\jac{M}{f(x,y)}_3\jac{1+2\w}{f(x,y)}_3=\jac{M}{f(x,y)}_3.
\]
Also,
\[
\jac{M}{a}_3^{-1} = \jac{\b_1by+3y_1\a(1+2\w)}{a}_3
=\jac{\b by+3y\a(1+2\w)}{a}_3,
\]
so
\[
\jac{M}{a}_3^{-1} =\jac{3A(1+2\w) +b(B/C)}{a}_3=
\jac{1+2\w)}{a}_3 R=R.
\]
Since $L=\jac{M}{af(x,y)}_3 \jac{M}{a}_3^{-1}$,
it remains to show that 
$\jac{M}{af(x,y)}_3=1$.
This will be done in two cases.

\noindent
{\bf Case 1:} $3 \mid Q$ or $3 \nmid (\b Q)$.

By Lemma 2.1, $3 \nmid \b$.
Since $(\b_1,3y_1 \a)=1$, we have $(N(M), \b_1)=1$ by (2.4).
In view of (2.2), we also have $(N(M),3y_1)=1$.
Cubic reciprocity yields
\begin{equation}\label{ eq. 2.6}
\jac{\b_1}{M}_3 = \jac{M}{\b_1}_
3=\jac{1+2\w}{\b_1}_3\jac{3y_1\a}{\b_1}_3=1.
\end{equation}
By (1.4) and (2.6),
\[
\jac{M}{af(x,y)}_3=\jac{af(x,y)}{M}_3 =
\jac{\b_1^2(ax+by)^2 - \b_1^2 y^2 d}{M}_3.
\]
By (2.4), 
\[
\b_1^2(ax+by)^2 \equiv -27 y_1^2\a^2 \pmod{N(M)},
\]
so since $\b_1 y=\b y_1$,
\[
\jac{M}{af(x,y)}_3=\jac{27y_1^2 \a^2 + \b^2 y_1^2 d}{M}_3
=\jac{27y_1^2 N(u)/Q^2}{M}_3=\jac{y_1}{M}_3^2\jac{Q}{M}_3.
\]

We next show that
\begin{equation}\label{ eq. 2.7}
\jac{Q}{M}_3 = 1.
\end{equation}

Let $q$ be a prime factor of $Q$.  If $q=3$, then
\[
\jac{q}{M}_3 = \jac{1+2\w}{M}_3^2 = \w^{4\a y_1}=1,
\]
since $3 \mid \a$ by Lemma 2.1.  If $q > 3$,  then
since $q \mid \a$ and $q \nmid \b$ by Lemma 2.1,
and since $q \nmid (ax+by)$ by (2.2), we have
\begin{equation}\label{eq. 2.8}
\jac{q}{M}_3 = \jac{M}{q}_3 =\jac{\b_1(ax+by)}{q}_3 = 1.
\end{equation}
Now suppose that $q=2$.  First assume that $D$ is even.
Then as in the proof of (2.5), 
we see that $\a$ is even and $\b(ax+by)$ is odd,
so that (2.8) holds for $q=2$.
Next assume that $D$ is odd.   Then as in the proof of (2.5),
$\a$ and $\b$ are odd while $(ax+by)$ and $y$ have opposite parity, 
so that either (2.8) holds for $q=2$ or else $(ax+by)$ is even and 
\[
\jac{2}{M}_3 = \jac{3y_1\a(1+2\w)}{2}_3 = 1.
\]
Thus $\jac{q}{M}_3=1$ for all primes $q$ dividing $Q$, which
completes the proof of (2.7).

To complete the proof for Case 1, it remains to prove that $\jac{y_1}{M}_3=1$.
Write $y_1=3^e y_0$, with $3 \nmid y_0$.  Then
\[
\jac{y_1}{M}_3 = \jac{3}{M}_3^e \jac{y_0}{M}_3
=\jac{1+2\w}{M}_3^{2e} \jac{M}{y_0}_3.
\]
Thus
\[
\jac{y_1}{M}_3 =\w^{4\a y_1 e} \jac{\b_1(ax+by)}{y_0}_3
=\w^{4\a y_1 e} =1,
\]
since $3 \mid y_1 e$.
This completes the proof for Case 1.

\noindent
{\bf Case 2:} $3 \nmid Q$ and $3 \mid \b$.

Since $3 \mid (B/C)$ and $3 \nmid Q$, we have $3 \nmid A$.
Recall that $u_2 = (A_2 + B_2 \sqrt{D})/T \in \S_D$,
where $A_2 = BD^2$, $B_2=AD$, and $T$ is the largest integer
cube dividing $(A_2, B_2)$.  Observe that $(3C,T)=1$,
since $3 \nmid A$ and every prime factor of $T$ must divide $D$.

Define $Q_2:=(A_2, B_2/C)/T$.  If $3 \mid D$, then $3 \mid Q_2$.
If $3 \nmid D$, then, since $3 \nmid A$, we have $3 \nmid (B_2Q_2)$.
Since (1.7) holds with $u_2$ in place of $u$,
Case 1 applies to show that $L(u_2) = R(u_2)$.   Then by
Lemmas 2.2--2.3,
\[
L(u) = L(u_2) = R(u_2) = R(u).
\]
This completes the proof for Case 2.

\end{proof}

\section{Proof of Theorem 1.1}

\begin{proof}
Let $u=A+B\sqrt{D} \in \S_D$ and consider a prime
\[
p=f(x,y) = ax^2 +2bxy+cy^2 \equiv 1 \pmod{3}
\]
with
$(ap,3yN(u)D)=1$, where the form $[a,2b,c]$ is in $H$.
In view of Theorem 2.4, (1.2), and (1.6), 
Theorem 1.1 is equivalent to
\begin{equation}\label{eq 3.1}
u  \mbox{ is a cubic residue} \pmod{p} \iff L(u)=1.
\end{equation}
By (2.1), 
\begin{equation}\label{eq 3.2}
L(u) = \jac{E}{p}_3, \quad \mbox{with  } E=3A+B\tau(1+2\w),
\end{equation}
where $\tau$ is a square root of $-3D$ mod $p$.
Since $p \equiv 1 \pmod{3}$, we can write $p=\p \bp$ for
Eisenstein primes $\p, \bp$.  By definition of the cubic residue symbol,
\[
\jac{E}{\bp}_3 \equiv E^{(p-1)/3}  \pmod{\bp}.
\]
Thus
\[
\jac{E}{\bp}_3^{-1} \equiv \bE^{(p-1)/3}  \pmod{\p},
\]
so that
\[
\jac{E}{\bp}_3 \equiv (1/\bE)^{(p-1)/3}    \pmod{\p}.
\]
Therefore
\[
L(u) = \jac{E}{p}_3 = \jac{E}{\p}_3\jac{E}{\bp}_3 
\equiv (E/\bE)^{(p-1)/3}  \pmod{\p}.
\]
By (3.2), $E \equiv 3A \pm 3B\sqrt{D} \pmod{\p}$, so that
\[
L(u) \equiv (u/u')^{\pm (p-1)/3} \pmod{\p},
\]
where $u' = A-B\sqrt{D}$.
Consequently,
\[
L(u)=1 \iff L(u) \equiv 1 \pmod{\p} \iff (u/u')^{(p-1)/3} 
\equiv 1 \pmod{\p}.
\]
Since $N(u)^{(p-1)/3} \equiv 1 \pmod{\p}$,
\[
L(u)=1 \iff (u^2)^{(p-1)/3} \equiv 1 \pmod{\p}
\iff u^{(p-1)/3} \equiv 1 \pmod{p}.
\]
Therefore $L(u)=1$ if and only if $u$ is a cubic residue mod $p$,
as desired.
\end{proof}

\section{A homomorphism on $H$}

For $u \in \S_D$, define a map $J: H \rightarrow \{1,\w,\w^2\}$
by $J(\C)=R(u,f)$, where $f=(a,2b,c)$ is a form in $\C$
with $(a,3N(u)D)=1$.  The map $J$ is well-defined on $H$ by the
same proof given in \cite{Sun2} except with Theorem 2.4 in place
of \cite[Theorem 2.1]{Sun2}. Moreover, $J$ is a group homomorphism.
To see this, let $\C_i = [a_i, 2b_i, c_i] \in H$ for  $i=1,2$,
where $(a_1a_2,3N(u)D)=1$ and $(a_1,a_2)=1$.
Then by Dirichlet composition \cite[Prop. 3.8]{Cox},
$\C_1\C_2 = \C_3$, where $\C_3=[a_3, 2b_3, c_3] \in H$ is given by
\begin{equation}\label{eq. 4.1}
a_3=a_1a_2, \quad b_3\equiv b_1 \pmod{a_1}, \quad
b_3 \equiv b_2 \pmod{a_2}.
\end{equation}
By (1.6) and (4.1), 
\[
J(\C_i) = \jac{9A -(B/C)b_3(1+2\w)}{a_i}_3, \quad i=1,2.
\]
Since $a_3=a_1a_2$, we see that $ J(\C_1)J(\C_2)=J(\C_3)$,
so $J$ is a homomorphism.

We close with the explanation promised in Example 1.2.
Let $H_2$ denote the subgroup of index 2 in $H$ consisting
of the classes that represent primes congruent to $1$ mod $3$.
(To see that $H_2$ is a group, suppose that
$\C_1, \C_2$ lie in $H_2$. Then $a_1\equiv a_2 \equiv 1 \pmod{3}$,
since the right member of (1.4) is congruent to $1$ mod $3$.
Thus $a_3=a_1a_2 \equiv 1 \pmod{3}$, so $\C_3 \in H_2$).
Let $G$ be the set of classes in $H_2$ that represent primes $p$
for which $u$ is a cubic residue mod $p$.  By Theorem 1.1, $G$ is the kernel
of the homomorphism $J: H_2 \rightarrow \{1,\w,\w^2\}$.
We proceed to show that this map is an epimorphism, i.e.,
\begin{equation}\label{eq. 4.2}
H_2/G \simeq \{1, \w, \w^2\}. 
\end{equation}
This will establish the promised result that
$G$ is a subgroup of index $6$ in $H$.

To prove (4.2), we will invoke the Cebotarev density theorem
\cite[Theorem 8.17]{Cox}.   
It suffices to show that
$u$ is a cubic residue mod $p$ for primes $p$ represented only by
one third of the classes in $H_2$,
or put another way, 
$u$ is a cubic residue mod $p$ for only ``a third"  of the 
primes $p$ satisfying $(-3/p) = (d/p) = (D/p)=1$.
Consider the number fields $K=\Q(\w,\sqrt{D})$ and $L= K(v)$,
where $v$ is the real cube root of $u$.  Recall that
$u$ is not a cube in $\OO$. Thus $u$ cannot be a cube
in $K$, so that  $|\Gal(L/K)|=3$.  
Let $W$ be the set of first degree prime ideals $\PP$ in $K$ which 
split completely in $L$.  
The Artin symbol $\left(\frac{L/K}{\PP}\right)$ is trivial for
$\PP \in W$, so 
by the Cebotarev density theorem,
$W$ has Dirichlet density 1/3.  
We can now relate $W$ to cubic residuacity, 
just as in the proof of \cite[Theorem 1.6]{ELSV},
to conclude that
$u$ is a cubic residue mod $p$ for primes $p$ represented only by
one third of the classes in $H_2$.

\section*{Acknowledgement}
The authors are grateful to Zhi-Hong Sun and Franz Lemmermeyer
for helpful suggestions.

\end{document}